\def\versiondate{27 June 2002}
\input math.macros
\input Ref.macros

\hyperstrue
\checkdefinedreferencetrue
\continuousfigurenumberingtrue
\theoremcountingtrue
\sectionnumberstrue
\forwardreferencetrue
\initialeqmacro

\def\path{{\cal P}}
\def\Pa{\gamma}
\def\F{{\cal F}}
\def\G{{\cal G}}
\def\wF{\widehat{\cal F}}
\def\wG{\widehat{\cal G}}
\def\B{{\cal B}}
\def\HH{{\cal H}}
\def\FF(#1, #2, #3, #4){{}_{#1}^{#2}\F_{#3}^{#4}}
\def\SS{{\cal S}}
\def\hit{{\ss hit}}
\def\LE{{\ss LE}}          
\def\eqD{{\buildrel {\cal D} \over =}}

\ifproofmode \relax \else\head{To appear in {\it Ann. Inst. H. Poincar\'e
Probab. Statist.}} {Version of \versiondate}\fi
\vglue20pt

\title{Markov Chain Intersections and the
 Loop-Erased Walk}

\author{Russell Lyons, Yuval Peres and Oded Schramm}

\abstract{Let $X$ and $Y$ be independent
transient Markov chains on the same state space that have
the same transition probabilities. Let $L$ denote the ``loop-erased
path'' obtained from the path of $X$ by erasing cycles when they are created.
We prove that if the paths of $X$ and $Y$ have
infinitely many intersections a.s., then $L$ and $Y$
also have infinitely many intersections a.s.}

\bottomII{Primary 
60J10. 
Secondary 
60G17
.} {Random walk, uniform
spanning forests.}
{Research partially
supported by the Institute for Advanced Studies, Jerusalem
(Lyons), NSF grant DMS--9803597 (Peres), and the Sam and Ayala
Zacks Professorial Chair (Schramm).}

\bsection {Introduction}{s.intro}

Erd\"os and Taylor (1960) proved that two independent simple
random walk paths in $\Z^d$ intersect infinitely often if $d \le
4$, but not if $d >4$.
Lawler (1991) proved that for $d=3,4$, there are still infinitely
many intersections even if one of the paths is replaced by
its loop-erasure, which is obtained by erasing cycles
in the order in which they are created.
Lawler's proof relied on
estimates that are only available in Euclidean space, so it
remained unclear how general this phenomenon is. Our main result,
\ref t.loop/ below, extends Lawler's result to general 
transient Markov chains.
Our primary motivation for studying intersections of  Markov chain
paths and loop-erased paths is the connection with uniform
spanning forests, which will be recalled in the next section.
The precise definition of loop-erasure is as follows.
We denote sequences by angle brackets $\Seq{\cdots}$ and sets by braces
$\{\cdots\}$.
Cardinality is denoted by bars $|\cdots|$.

\medbreak\definition Let $\SS$ be a countable set, and
consider a (finite or infinite) sequence $\Pa=\Seq{v_0,v_1, v_2,
\ldots}$ in $\SS$, such that each $v\in\SS$ occurs only finitely often
in $\Pa$. The {\bf loop-erasure} of $\Pa$, denoted 
$\LE(\Pa) = \Seq{u_0, u_1, u_2, \ldots}$, is constructed
inductively as follows.
Set $u_0:=v_0$. If $u_j$ has been chosen, then let $k$ be the last
index such that $v_k=u_j$. If $v_k$ is the last element of $\Pa$,
then let $\LE(\Pa):=(u_1,\dots,u_j)$;
 otherwise,  set $u_{j+1}:=v_{k+1}$ and continue.
\medbreak

We shall let $\LE(\Pa)$ stand for the sequence as well as the set.
In this notation, our main result can be stated as follows:

\procl t.loop Let $\Seq{X_m}$ and $\Seq{Y_n}$ be independent
transient Markov chains on the same state space $\SS$ that have
the same
transition probabilities, but
possibly different initial states. Then on the event that $|\{X_m\} \cap
\{Y_n\}| = \infty$, almost
surely $|\LE{\Seq{X_m}} \cap \{Y_n\}| = \infty$.
\endprocl

The key to proving this theorem is the following quantitative lemma.

\procl l.IntLoop Let $\Seq{X_m}$ and $\Seq{Y_n}$ be independent transient
Markov chains
on the same state space $\SS$ that have the same transition probabilities,
but possibly different initial states.
Then 
$$
\P[\LE\Seq{X_m} \cap \{ Y_n \} \ne \emptyset ]
\ge
2^{-8}\, \P[\{X_m\} \cap \{ Y_n \} \ne \emptyset ]
\,.
$$
\endprocl

\procl r.kill
\ref l.IntLoop/ also applies to Markov chains that are
killed when they exit a certain set, or killed at an exponential time;
in this form, it contains new information even when the underlying
chain is simple random walk in $\Z^d$ for $d \ge 3$.
\endprocl

\ref t.loop/ reduces the question of the intersection of a Markov chain and
the loop-erasure of an independent copy to the simpler problem of the
intersection of two independent copies of the chain. However,
for a general Markov chain, it may still be difficult to find the probability
that
two independent sample paths of the chain have infinitely many intersections.
Consider the Green function
$$
G(x,y):=\sum_{n=0}^\infty \P_x[X_n=y] \,.
$$
Because of the Borel-Cantelli lemma, if the
independent sample paths $X$ and $Y$ satisfy
\begineqalno
\sum_{m,n=0}^\infty \P_{x,x}[X_m=Y_n]
 &= \sum_{z \in
\SS}\sum_{m,n=0}^\infty \P_{x}[X_m=z]\P_{x}[Y_n=z] \cr &= \sum_{z
\in \SS} G(x,z)^2 <\infty \, ,  \label e.Green    \cr
\endeqalno
then the number of intersections of $X$ and $Y$ is a.s.\ finite. In
general, the converse need not hold; see \ref x.counter/.
Nevertheless, if the transition probabilities are invariant under
a transitive group of permutations of the state space, then the
converse to \ref e.Green/ is valid.

\procl t.trans Let $p(\cdot,\cdot)$ be a
transition kernel on a countable state space $\SS$. Suppose that
$\Pi$ is a group of permutations of $\SS$ that acts
transitively ({\it i.e.}, with a single orbit) and satisfies
$p(\pi x, \pi y)=p(x,y)$ for all $\pi \in \Pi$ and
$x,y \in \SS$. Suppose that
 $$ \sum_{z \in \SS} G(o,z)^2 =\infty \, ,
\label e.Green2 $$
where $o$ is a fixed element of $\SS$.
Then two independent chains $X$ and $Y$ with
transition probabilities $p(\cdot,\cdot)$ and initial state $o \in
\SS$ have infinitely many intersections a.s. Moreover,
 $Y_n$ is in
$\LE(\Seq{X_m}_{m \ge 0})$ for infinitely many $n$ a.s.
\endprocl

For simple random walk on a vertex-transitive graph,
 \ref e.Green2/ holds if and only if the graph has polynomial volume
 growth of degree at most $4$; see \ref c.transgraph/.

Next, we consider triple intersections. 
It is well known that three independent
simple random walks in $\Z^3$ have infinitely many
mutual intersections a.s.\ (see Lawler (1991), Section 4.5).
To illustrate the versatility of \ref t.loop/,
we offer the following refinement.

\procl c.three
 Let $X=\Seq{X_m}, \; Y=\Seq{Y_n}$  and $Z=\Seq{Z_k}$ be independent
simple random walks in the lattice $\Z^3$. Denote by $L_Z(X)$
the ``partially loop-erased''  path,
 obtained from $X$ by erasing any cycle that starts (and ends)
at a node in $Z$, where the erasure is made when the cycle is created.
Then the (setwise) triple intersection $L_Z(X) \cap Y \cap Z$ is a.s.\ infinite.
\endprocl

 See \ref c.three2/ in
\ref s.conc/ for an extension and the (very short) proof.
We note that $L_Z(X)$ cannot be replaced by $\LE(X)$ in this
corollary; this follows from Lawler (1991), Section 7.5.

\medskip

In the next section, we recall the connection to spanning forests
and  sketch a heuristic argument for \ref t.loop/.
In \ref s.second/, we
discuss the reverse second moment method for two Markov chains,
following Salisbury (1996). \ref l.IntLoop/ and
\ref t.loop/ are proved in \ref s.loop/
by combining ideas from the two preceding sections.
\ref s.trans/ contains a proof of \ref t.trans/ on
transitive chains.
 Concluding remarks and questions are in \ref s.conc/.

\bsection{Spanning Forests and Heuristics}{s.heur}

Loop-erased random walks and uniform spanning trees are intimately related.
Let $G$ be a finite graph with $x, y$ two vertices of $G$.
Let $L$ be the loop-erasure of the random walk path started at $x$ and
stopped when it reaches $y$.
On the other hand, let $T$ be a spanning tree of $G$ chosen uniformly
and let $L_T$ be
the shortest path in $T$ that connects $x$ and $y$.
Pemantle (1991) showed that $L_T$ has the same distribution as $L$.
Given that $L_T = \ell$ for some simple path $\ell$, the remainder of
$T$ has the uniform distribution among spanning trees of the graph obtained
from $G$ by contracting $\ell$.
Therefore,
it follows immediately from $L_T \eqD L$ 
that a uniform spanning tree can also be chosen as follows.
Pick any vertex $x_0$ of $G$. Let $L_0$ be loop-erased random walk from any
vertex $x_1$ to $x_0$.
Pick any vertex $x_2$ and let $L_1$ be loop-erased random walk from $x_2$ to
$L_0$.
Pick any vertex $x_3$ and let $L_1$ be loop-erased random walk from $x_2$ to
$L_0 \cup L_1$.
Continue until a spanning tree $T := L_0 \cup L_1 \cup \cdots$ is created.
Then $T$ has the uniform distribution.

This is known as Wilson's algorithm for generating a uniform spanning tree.
Wilson (1996) showed that an analogous algorithm exists corresponding to any
Markov chain, not merely to reversible Markov chains.

Next, we discuss the analogous object on infinite graphs.
The {\bf wired uniform spanning forest} (WUSF) in an infinite graph $G$ may be
defined as a weak limit of uniform random spanning trees in an
exhaustion of $G$ by  finite subgraphs $G_n$, with the boundary of
$G_n$ identified to a single point (``wired''). The resulting
measure on spanning forests does not depend on the exhaustion. The
WUSF was implicit in Pemantle (1991)
 and was made explicit by H\"aggstr\"om (1995); see Benjamini,
Lyons, Peres, and Schramm (2001), denoted BLPS (2001) below, for
details. The connection of the WUSF to loop-erased walks was
discovered by Pemantle (1991):

\procl p.onetree
Let $G$ be a locally-finite connected graph.
The wired uniform spanning forest (WUSF) is a single tree a.s.\
iff from every (or some) vertex, simple random walk and an
independent loop-erased random walk intersect infinitely often
a.s.  Moreover, the probability that $u$ and $v$ belong to the
same WUSF component equals the probability that a simple random
walk path from $u$ intersects an independent loop-erased walk from
$v$.
\endprocl

 Just as the relation between spanning trees  and
loop-erased walks in finite graphs was clarified by the algorithm
of Wilson (1996) for generating uniform spanning trees, this
algorithm was extended to infinite graphs in BLPS (2001) to generate the WUSF.
With this extended algorithm, \ref p.onetree/ becomes obvious.
This proposition illustrates why \ref t.loop/ is useful in the study of
the WUSF.

 We now sketch a heuristic argument for  \ref t.loop/.
On the event
that $X_m=Y_n$, the continuation paths $X':=\Seq{X_j}_{j \ge m}$ and
$Y':=\Seq{Y_k}_{k \ge n}$ have the same distribution, whence the
chance is at least $1/2$ that $Y'$ intersects $L:=\LE\Seq{X_0,\ldots ,X_m}$ at
an earlier (in the clock of $L$) point than $X'$.
On this event, the earliest intersection point of $Y'$ and $L$
will remain in $\LE\Seq{X_j}_{j \ge 0} \bigcap \Seq{Y_k}_{k \ge
0}$. The difficulty in making this heuristic precise lies in
selecting a pair $(m,n)$ such that $X_m=Y_n$, given that such pairs
exist. The natural rules for selecting such a pair ({\it e.g.},
lexicographic ordering) affect the law of  at least one of the
continuation paths, and invalidate the argument above;
\msnote{appendix}
R. Pemantle (private communication, 1996) showed that this holds
for {\it all\/} selection rules. Our
solution to this difficulty is based on applying a second moment
argument to a weighted count of intersections.

\bsection{The Second Moment Method and a Converse}{s.second}

\procl t.seconv Let $X$ and $Y$ be two independent Markov chains
on the same countable state space $\SS$, with initial states $x_0$
and $y_0$, respectively.  Let
$$
A \subset \N \times \SS \times
\N \times \SS\,,
$$
and denote by  $\hit(A)$  the event that $(m, X_{m}, n, Y_{n}) \in A$
for some $m,n \in \N$.
Given any weight function $ w:  \SS\to [0,\infty)$
that vanishes outside of $A$,
consider the
random variable $$ S_w:=\sum_{m,n=0}^\infty w(m,X_m,n,Y_n) \,.$$
If $\P[\hit (A)]>0$, then there exists
such a $w$ satisfying $0<\E [S_w^2] < \infty$ and $$
\P[\hit(A)] \le  64 {(\E S_w)^2 \over \E[S_w^2]} \,.
\label e.ii
$$
\endprocl

Note that this provides a converse estimate to that provided by the 
Cauchy-Schwarz inequality (often referred to as ``the second-moment
method"):
If $0<\E [S_w^2] < \infty$, then
$$ \P[\hit(A)] \ge \P[S_w>0] \ge {(\E S_w)^2 \over \E[S_w^2]} \,.
\label e.i
$$

\ref  t.seconv/  is essentially contained in Theorem 2 of Salisbury (1996),
which is, in turn, based on the ideas in the path-breaking paper of
Fitzsimmons and Salisbury (1989).
We include the proof of this theorem,
 since our focus on time-space chains allows
us to avoid the subtle time-reversal argument in Salisbury's paper.
The ratio $64$ between the upper and lower bounds in \ref e.ii/
and \ref e.i/, respectively, 
improves on the ratio $1024$ obtained in Salisbury
(1996), but we suspect that it is still not optimal.
We remark that the lower bound for the hitting probability
which is stated but not proved in Corollary 1 of Salisbury (1996)
(i.e., the left-hand inequality in the last line of the statement)
is incorrect as stated, but we shall not need that inequality.

We start with two known lemmas.
\procl l.3fields Let $(\Omega,\B, \P)$ be a probability space.
Suppose that $\F,\G,\HH$ are sub-$\sigma$-fields of $\B$ such that
$\G \subset \F \cap \HH$ and
$\F, \HH$ are conditionally independent given $\G$. Then
\beginitems
\itemrm{(i)} for any $f \in L^1(\F):= L^1(\Omega, \F, \P)$,
$$
\E[f \mid \HH]=\E[f \mid \G] \,;
\label e.HG
$$
\itemrm{(ii)}
for any  $\varphi \in L^1(\B):=L^1(\Omega, \B, \P)$,
$$ \Ebig{ \E[\varphi \mid \F] \bigm|
\HH } = \E[\varphi \mid \G] \,.
\label e.varphi
$$
\enditems
\endprocl

\proof \item{(i)} For any $h \in L^\infty(\HH)$
and $g \in L^1(\G)$, we have
$\int hg \, d\P= \int \E[h \mid \G]g \, d\P$
 by definition of conditional expectation. In particular,
$$
\int h \E[f \mid \G] \, d\P = \int \E[h \mid \G] \E[f \mid \G]\, d\P
 =\int \E[hf \mid \G] \, d\P
=\int \! h f \, d\P \, ,
\label e.iden
$$
where the second equality follows from the conditional
independence assumption. Since $\E[f \mid \G]$ is $\HH$-measurable,
the identity \ref e.iden/ for all  $h \in L^\infty(\HH)$ implies
that \ref e.HG/ holds.
\item{(ii)} Write $f:= \E[\varphi \mid \F]$. Since
$\E[f \mid \G] =\E[\varphi \mid \G]$, \ref e.varphi/ follows from
\ref e.HG/.
\Qed

We shall use the following inequality from Burkholder,
Davis and Gundy (1972).

 \procl l.bdg Suppose that $\Seq{\F_m}$ is an increasing
or decreasing sequence of $\sigma$-fields and that $\Seq{\varphi(m)}$ is a
sequence of nonnegative random variables. Then
$$
\Eleft{ \biggl( \sum_m \E[\varphi(m) \mid \F_m] \biggr)^2 }
\le 4 \, \Eleft{ \biggl( \sum_m \varphi(m) \biggr)^2 }\,.
$$
\endprocl

\proofof t.seconv
 Consider the $\sigma$-fields
$$\eqalign{
\F_m:=\sigma(X_1,\ldots , X_m) \, , & \qquad
\wF_m:=\sigma(X_m, X_{m+1}, \ldots ) \, , \cr
\G_n:=\sigma(Y_1,\ldots , Y_n) \, , & \qquad \wG_n:=\sigma(Y_n,
Y_{n+1}, \ldots ) \, . \cr}
\label e.sfield
$$
Abbreviate $\F:=\wF_0$ and $\G=\wG_0$.

We begin with the lexicographically minimal stopping time 
$(\tau,\lambda)$ defined as follows:
\vskip\baselineskip

\vbox{
\begingroup\narrower\narrower
\noindent
if there exist $m,n \ge 0$ such that $(m,X_m,n,Y_n) \in A$,
then let
\begineqalno
\tau&:=\min\{m \st \exists n \; \; (m,X_m,n,Y_n) \in A\}
\,,\cr
\lambda&:=\min\{n \st  (\tau,X_\tau,n,Y_n) \in A \} \,;
\endeqalno
otherwise, set $\tau:=\lambda:=\infty$.

\endgroup
}

\noindent
Consider
$$
\psi(m,n):=\II{\tau=m, \lambda=n} \,.
$$
Since $\psi(m,n)$ is $\F_m \vee \G$-measurable and
$\F_m \vee \G$ is conditionally  independent of
$\wF_m \vee \G$  given $\sigma(X_m) \vee \G$,
\ref l.3fields/(i) implies that
$$
\psi_1(m,n):=\E[\psi(m,n) \mid \wF_m \vee \G]
=\E[\psi(m,n) \mid \sigma(X_m)\vee \G] \,.
$$
Let $\psi_2(m,n):=\E[\psi_1(m,n) \mid \F \vee \G_n]$.
Two applications of \ref l.3fields/(ii) yield that
$$
\psi_3(m,n):=\E[\psi_2(m,n) \mid \F \vee \wG_n]
= \E[\psi_1(m,n) \mid \F \vee \sigma(Y_n)]=
\E[\psi(m,n) \mid \sigma(X_m, Y_n)] \,.
$$
Thus
$$
\psi_3(m,n)=\P[\tau=n, \lambda=m \mid X_m, Y_n] = 
\sum_{x,y \in \SS} w(m,x,n,y)\II{X_m=x, Y_n=y}\,,
\label e.psi3
$$
where
 $$
w(m,x,n,y):=\P[\tau=m, \lambda=n \mid X_m=x, Y_n=y] \,.
$$
Applying \ref l.bdg/ to the random variables
$
\varphi(m):=\sum_{n \ge 0} \psi(m,n) 
$
and the $\sigma$-fields $\wF_m \vee \G$, we obtain
$$
\Eleft{ \biggl( \sum_{m,n} \psi_1(m,n) \biggr)^2 }
\le 4 \, \Eleft{ \biggl( \sum_{m,n} \psi(m,n)  \biggr)^2 }
= 4 \, \P[\hit (A)]
\,,
$$
since  $\sum_{m,n} \psi(m,n)=\I{\hit (A)}$.
By \ref e.psi3/, we have
$$
S_w:=\sum_{m,n} w(m,X_m,n,Y_n)= \sum_{m,n} \psi_3(m,n) \,.
$$
Two applications of \ref l.bdg/,
first with the variables
$\varphi_2(n):=\sum_m \psi_2(m,n)$
and the $\sigma$-fields $\F \vee \wG_n$,
then with the variables
$\varphi_1(n):=\sum_m \psi_1(m,n)$
and the $\sigma$-fields $\F \vee \G_n$, yield
$$
\E[S_w^2] \le
4\Eleft{ \biggl( \sum_{m,n} \psi_2(m,n) \biggr)^2 } \le
16 \, \Eleft{ \biggl( \sum_{m,n} \psi_1(m,n) \biggr)^2 }
 \le 64 \, \P[\hit (A)] \,.
$$
Since $\E[S_w]=  \P[\hit (A)]$,  the previous inequality
is equivalent to \ref e.ii/.
\Qed

\bsection{Intersecting the Loop-Erasure}{s.loop}

We shall prove the following extension of \ref l.IntLoop/.
\procl l.IntLoop2
 Let $\Seq{X_m}_{m \ge 0}$ and $\Seq{Y_n}_{n \ge 0}$
be independent transient Markov chains
on $\SS$ that have the same transition probabilities,
but possibly different initial states $x_0$ and $y_0$.
Given $k \ge 0$, fix $\Seq{x_j}_{j=-k}^{-1}$ in $\SS$ and set
$X_j:=x_j$ for $-k \le j \le -1$.
Then the probability that the loop-erasure of
$\Seq{X_m}_{m \ge -k}$ intersects $\{Y_n\}_{n \ge 0}$
is at least
$2^{-8}\, \P[\texists {m \ge 0}
\texists {n \ge 0} X_m=Y_n ]$.
\endprocl
\proof
For $A:=\{(m,x,n,x) \st m,n \ge 0, x \in \SS\}$, choose a weight
function $w:A \to [0,\infty)$ as in \ref t.seconv/, which defines the sum
$S_w$.
Denote
$$
L^m:=\Seq{L^m_j}_{j=0}^{J(m)}:=\LE\Seq{X_{-k},X_{1-k},\ldots,X_m} \,.
$$
On the event $\{X_m=Y_n\}$, define
\begineqalno
&j(m,n):=\min\bigl\{j \ge -k \st
 L^m_j \in  \{X_m,X_{m+1},X_{m+2}, \ldots\}
\bigr\},  \label e.defj \cr
&i(m,n):=\min\bigl\{i  \ge -k
\st L^m_i \in  \{Y_n,Y_{n+1},Y_{n+2}, \ldots\}
\bigr\} \, .
\label e.defi
\endeqalno
Note that the sets on the right-hand sides of \ref e.defj/
and \ref e.defi/ both  contain $J(m)$ if $X_m=Y_n$.
Define $j(m,n):=i(m,n):=0$ on the event $\{X_m \neq Y_n\}$.
Let $\chi(m,n):=1$ if $i(m,n) \le j(m,n)$, and $\chi(m,n):=0$ otherwise.
Given $\{X_m=Y_n=x\}$, the continuations
$\Seq{X_m,X_{m+1},X_{m+2}, \ldots}$ and $\Seq{Y_n,Y_{n+1},Y_{n+2},
\ldots}$ are exchangeable with each other, so for every $x \in \SS$,
$$
\Ebig{\chi(m,n) \mid X_m=Y_n=x} =
\Pbig{i(m,n) \le j(m,n) \mid X_m=Y_n=x} \ge {1 \over 2} \,.
\label e.key1
$$
Observe that if $X_m=Y_n$ and $i(m,n) \le j(m,n)$, then
$L^m_{i(m,n)}$ is in $\LE\Seq{X_r}_{r=-k}^\infty
 \bigcap \{Y_\ell\}_{\ell=0}^\infty $.

Consider the random variable
$$
\Upsilon_w:=\sum_{m=0}^\infty \sum_{n=0}^\infty
w(m,X_m,n,Y_n)\,\chi(m,n) \,.
$$
Obviously $\Upsilon_w \le S_w$ everywhere. On the other hand,
by conditioning on $X_m, Y_n$ and applying \ref e.key1/, we see that
$$
\E [\Upsilon_w] =\sum_{m=0}^\infty \sum_{n=0}^\infty
\EBig{  w(m,X_m,n,Y_n)\, \Ebig{\chi(m,n) \mid X_m,Y_n} }
\ge {1 \over 2} \E[S_w] \,.
$$
By our choice of $w$ and \ref t.seconv/, we conclude that
$$
P[\Upsilon_w>0] \ge {(\E \Upsilon_w)^2 \over \E[\Upsilon_w^2]}
\ge {(\E S_w)^2 \over 4 \E[S_w^2]}
\ge  {1 \over 256} \P[\hit(A)] \,.
$$
The observation following \ref e.key1/ and the definition of $A$
conclude the proof.
\Qed

The next corollary follows immediately from \ref l.IntLoop2/
and the Markov property of $\Seq{X_m}$ at a fixed
time $k$.
\procl c.cond
 Let $\Seq{X_m}_{m \ge 0}$ and $\Seq{Y_n}_{n \ge 0}$
be independent transient Markov chains
on $\SS$ that have the same transition probabilities,
but possibly different initial states $x_0$ and $y_0$.
Suppose that the event $B:=\{X_1=x_1,\ldots,X_k=x_k\}$
has $\P[B]>0$. Then
$$
\P_{x_0,y_0}
\biggl[\LE\Seq{X_m}_{m \ge 0} \cap \{Y_n\} \neq \emptyset \, \Big|
\, B \biggr] \ge
{1 \over 256} \P_{x_0,y_0}\biggl[\texists {m \ge k}
\texists {n \ge 0} X_m=Y_n \, \Big|\,  B \biggr]\,.
\label e.cond
$$
\endprocl

\proofof t.loop
Denote by $\Lambda$ the event that $Y_n$
is in  $\{X_m\}_{m \ge 0}$ for infinitely many $n$.
Let $\Gamma_\ell$ denote the event that
$Y_n \notin \LE(\Seq{X_m}_{m \ge 0})$ for all $n \ge \ell$,
and define $\Gamma:=\cup_{\ell \ge 0} \Gamma_\ell$.
We must show that  $\P_{x_0,y_0}[\Lambda \cap \Gamma]=0$
for any choice of initial states $x_0, y_0$.
Suppose, on the contrary, that
 $$
\exists x_0, y_0 \quad
\P_{x_0,y_0}[\Lambda \cap \Gamma]>0 \,.
\label e.contrary
$$
Then for $\ell$ sufficiently large,
$\P_{x_0,y_0}[\Lambda \cap \Gamma_\ell]>0$.
 By L\'evy's zero-one law,
for any $\epsilon>0$, if $k,r$ are large enough, then
 there exist $x_1,\ldots,x_k, \;  y_1,\ldots, y_r$
such that the events
$$
B:=\{X_1=x_1, \ldots, X_k=x_k\}, \quad B':=\{Y_1=y_1,\ldots,Y_r=y_r \}
$$
satisfy $\P_{x_0,y_0}[B \cap B']>0$ and
$\P_{x_0,y_0}[\Lambda \cap \Gamma_\ell \mid B \cap B']>1-\epsilon$.
We fix such events $B,B'$ with $r>\ell$.
Starting the chain $Y$ at $y_r$ instead of $y_0$
and using the Markov property,
we infer that
$$
\P_{x_0,y_r}[\Lambda \cap \Gamma_0 \mid B]> 1-\epsilon \, .
\label e.gaga
$$
However, \ref c.cond/ implies that
$$
\P_{x_0,y_r}[\Gamma_0^c \mid B] \ge {1 \over 256}
\P_{x_0,y_r}[\Lambda \mid B] >  {1-\epsilon \over 256}    \, .
$$
Adding the preceding two inequalities,
we obtain
$$
1 \ge \P_{x_0,y_r}[\Lambda \cap \Gamma_0 \mid B]+
 \P_{x_0,y_r}[\Gamma_0^c \mid B] >{257(1-\epsilon) \over 256} \,.
$$
Taking $\epsilon<1/257$ yields a contradiction
to the assumption \ref e.contrary/
and completes the proof.
\Qed

\bsection{Transitive Markov Chains}{s.trans}

One ingredient of the proof of \ref t.trans/ will be the following
lemma, which does not require transitivity.

\procl l.slide Let $X$ and $Y$ be independent transient
 Markov chains on the same state space that
have the same
transition probabilities.
 Denote by $\Lambda$ the event that
the paths of $X$ and $Y$ intersect infinitely often,
and let $u(x,y):=\P_{x,y}[\Lambda]$, where the subscripts indicate the
initial states of $X$ and $Y$ respectively.
Then $u(x,x) \ge 2u(x,y)-1$ for all $x, y$.
\endprocl
\proof  Since $u(x,\cdot)$ is harmonic, the sequence
$\Seq{u(x,Y_n)}_{n \ge 0}$
is a bounded martingale. Therefore
\begineqalno
u(x,y)-u(x,x)
&=\lim_{n \to \infty} \E_y[u(x,Y_n)]-\lim_{m \to \infty} \E_x[u(x,X_m)]
\cr
&=\E_{x,y}\Bigl[\lim_{n \to \infty} u(x,Y_n)-\lim_{m \to \infty}
u(x,X_m)\Bigr] 
\,.\label e.limm \cr
\endeqalno
On the event $\Lambda$, the  two
limits in \ref e.limm/
coincide; therefore,
$$
u(x,y)-u(x,x) \le \P_{x,y}[\Lambda^c]= 1-u(x,y) \,.
$$
This is equivalent
to the assertion of the lemma. \Qed

\proofof t.trans
Here both Markov chains $X,Y$ are started
at $o$, so  we write $\P$ rather than $\P_o$, etc.
 Denote $G_n(o,x):=\sum_{k=0}^n \P[X_k=x]$. By
transitivity,  $$\sum_{w \in \SS} G_n(z,w)^2 =\sum_{w \in \SS}
G_n(o,w)^2 \label e.allz $$ for all $z \in \SS$. Let
$I_n:=\sum_{k=0}^n \sum_{m=0}^n \II{X_k=Y_m}$ be the number of
intersections of $X$ and $Y$ by time $n$. Then
\begineqalno
 \E[I_n]
 &= \sum_{z \in \SS} \sum_{k=0}^n \sum_{m=0}^n \P[X_k=z=Y_m] \cr
 &= \sum_{z \in \SS} \sum_{k=0}^n  \P[X_k=z]
 \cdot \sum_{m=0}^n \P[Y_m=z] \cr
 &=\sum_{z \in \SS} G_n(o,z)^2 \,.
  \label e.first \cr
\endeqalno
To estimate the second moment of $I_n$, observe that
\begineqalno
  \sum_{k,i=0}^n \P[X_k=z, X_i=w]
 &=  \sum_{k=0}^n  \sum_{i=k}^n \P[X_k=z] \P[ X_i=w \mid X_k=z]
 \cr
 &\quad + \sum_{i=0}^n  \sum_{k=i+1}^n \P[X_i=w] \P[ X_k=z \mid X_i=w]
 \cr
 & \le G_n(o,z) \,G_n(z,w) + G_n(o,w) \,G_n(w,z) \,.  \cr
\endeqalno
Therefore
\begineqalno
 \E[I_n^2]
 &= \sum_{z,w \in \SS} \, \, \sum_{k,m=0}^n \, \, \sum_{i,j=0}^n \P[X_k=z=Y_m, X_i=w=Y_j ] \cr
 &= \sum_{z,w \in \SS} \, \, \sum_{k,i=0}^n  \P[X_k=z, X_i=w]
 \cdot \sum_{m,j=0}^n \P[Y_m=z, Y_j=w] \cr
 & \le \sum_{z,w \in \SS} [G_n(o,z) \,G_n(z,w) + G_n(o,w) \,G_n(w,z)]^2 \cr
 & \le \sum_{z,w \in \SS} 2 [G_n(o,z)^2 \,G_n(z,w)^2 + G_n(o,w)^2 \,G_n(w,z)^2] 
 \cr&= 4\sum_{z,w \in \SS} G_n(o,z)^2 \,G_n(z,w)^2 \,.
  \label e.second \cr
\endeqalno
Summing first over $w$ and using \ref e.allz/, then \ref e.first/,
we deduce that
 $$\E[I_n^2] \le 4 \Big(\sum_{z \in \SS} G_n(o,z)^2\Big)^2 =4\,\E[I_n]^2
 \,.
\label e.sec2
$$
By a consequence of the Cauchy-Schwarz inequality (see, e.g., Kahane
(1985), p.\ 8),
$$\P\Big[I_n \ge \epsilon \E[I_n]\Big] \ge
(1-\epsilon)^2 {\E[I_n]^2 \over \E[I_n^2]} \ge {(1-\epsilon)^2
\over 4} \, .
\label e.kahane
$$
As in \ref l.slide/, denote by $\Lambda$ the event that
the path-sets $X$ and $Y$ have infinitely many
intersections, and let $u(x,y):=\P_{x,y}[\Lambda]$.
Define $\F_n, \G_m$ as in \ref e.sfield/.
Apply the transience of $X$ and $Y$ and the Markov property
to obtain
$$
\P[\Lambda \mid \F_n \vee \G_n]= \P[\Lambda \mid X_n, Y_n]=
u(X_n,Y_n) \,.
$$
Therefore,
by L\'evy's zero-one law, $\lim_{n \to \infty} u(X_n,Y_n)=\I{\Lambda}$
a.s.

By \ref e.kahane/ and the hypothesis \ref e.Green2/, 
$ \P[\Lambda]=\P[\lim I_n =\infty] \ge 1/4$.
On the event $\Lambda$, we have  by \ref l.slide/ that
$$
\lim_{n \to \infty} u(X_n,X_n) \ge 2\lim_{n \to \infty} u(X_n,Y_n) -1
=1 \,,
$$
whence $u(o,o)= 1$ by transitivity. The assertion
concerning loop-erased walk now follows from \ref t.loop/.
\Qed
\procl r.gallr
The calculation leading to \ref e.sec2/ follows
Le Gall and Rosen (1991), Lemma~3.1. More generally, their argument
gives $\E[I_n^k] \le (k!)^2 (\E I_n)^k$ for every $k \ge 1$.
\endprocl

\procl c.transgraph
Let $\Delta$ be an infinite, locally finite,
vertex-transitive graph. Denote by $V_n$ the number 
of vertices in $\Delta$ at distance at most $n$ from a fixed vertex $o$.
\beginitems
\itemrm{(i)} If $\sup_n V_n/n^4= \infty$, then
two independent sample paths of simple random walk in 
$\Delta$ have finitely many intersections
a.s.
\itemrm{(ii)} Conversely, if $\sup_n V_n/n^4 <\infty$, then
two independent sample paths of simple random walk in 
$\Delta$ intersect infinitely often a.s.
\endprocl
\enditems
\proof
For independent simple random walks, reversibility and regularity of $\Delta$
imply that
$$
\sum_{m=0}^\infty \sum_{n=0}^\infty  \P_{x,x} [X_{m}=Y_n] =
\sum_{m=0}^\infty \sum_{n=0}^\infty  \P_x[X_{m+n}=x]= 
\sum_{n=0}^\infty (n+1) \P_x[X_{n}=x] \,.
\label e.heat
$$
\item{(i)} 
The assumption  that $\sup_n V_n/n^4= \infty$ 
implies that $V_n \ge cn^5$ for some $c>0$
and all $n$: see Theorem 5.11 in Woess (2000).
Corollary 14.5 in the same reference yields
$\P_x[X_n=x] \le Cn^{-5/2}$. Thus the sum in \ref e.heat/
converges.
\item{(ii)} Combining the results (14.5), (14.12)
and (14.19) in Woess (2000), we infer that the assumption
$V_n=O(n^4)$ implies that $P_x[X_{2n}=x] \ge cn^{-2}$
for some $c>0$ and all $n \ge 1$. Thus the series \ref e.heat/
diverges, so the assertion follows from \ref e.Green/
and \ref t.trans/.
\Qed

\bsection {Concluding Remarks}{s.conc}

The following example shows that the
invariance assumption in \ref t.trans/ cannot be omitted.
\procl x.counter
Consider the graph $H$, defined as the
union of copies of $\Z^5$ and $\Z$, joined at a single common vertex $o$.
The Green function  for simple random walk in  $H$ satisfies
$$ G(o,z) ={\deg z \over \deg{o}} G(z,o)= {2 \over 12} G(o,o) 
$$
provided $z \neq o$ is in the copy of $\Z$.
In particular, $\sum_z  G(o,z)^2 =\infty$.
However, two independent simple random walks on  $H$ will
have finitely many intersections a.s.
\endprocl

We continue by proving an extension of \ref c.three/.
\procl c.three2
 Let $X=\Seq{X_m}$ and  $Y=\Seq{Y_n}$  be independent
transient Markov chains on a state space $\SS$ and that have
the same 
transition probabilities.
Let $Z$ be a subset of
$\SS$ such that $Z$ is a.s.\ hit infinitely often by $X$ (and so by $Y$).
 Denote by $L_Z(X)$
the sequence obtained from $X$ by erasing any cycle that starts (and ends)
at a state in $Z$, where the erasure is made when the cycle is created.
Then on the event that  $X \cap Y \cap Z$ is infinite, almost
surely $L_Z(X) \cap Y \cap Z$ is also infinite.
\endprocl
\proof
Let $m(0)=0$ and $m(j+1):= \min\{k > m(j): X_k \in Z\}$
for all $j \ge 0$. Then $X^Z=\Seq{X^Z_j}:= \Seq{X_{m(j)}}$ is a Markov chain
(``the chain $X$ induced on $Z$'').
Similarly, let $Y^Z=\Seq{Y^Z_i}= \Seq{Y_{n(i})}$
denote the chain $Y$ induced on $Z$. Since $ \LE(X^Z) =L_Z(X) \cap Z$,
and $Y^Z=Y \cap Z$ as sets of vertices,
the assertion follows by applying \ref t.loop/
to the chains $X^Z$ and $Y^Z$ on the state space $Z$.
\Qed

  A natural question suggested by \ref t.loop/ is the following.
\procl q.looploop
Let $X=\Seq{X_m}$ and $Y=\Seq{Y_n}$ be independent transient
Markov chains on a state space $\SS$ and that have the same
transition probabilities. Suppose that
$|\{X_m\} \cap \{Y_n\}| = \infty$ a.s.
Must $|\LE{\Seq{X_m}} \cap \LE\Seq{Y_n}| = \infty$ a.s?
\endprocl
This question is open even for simple
random walk in $\Z^3$. For simple random walk in $\Z^4$, an
affirmative answer was given by Lawler (1998).

Our final question arose from an attempt to compare the stationary
$\sigma$-fields defined by a Markov chain and by its loop-erasure.
\procl q.tail
Let $\Seq{X_m}$ be a transient
Markov chain.
Consider $\Seq{L_j}=\LE\Seq{X_m}_{m \ge 0}$ and
 $\Seq{L_j^*}=\LE\Seq{X_m}_{m \ge 1}$. Does there a.s.\ exist
some integer $k$ such that $L_j^*=L_{j+k}$ for all large $j$?
\endprocl
The answer is certainly positive if  $\Seq{X_m}$
has  infinitely many ``cutpoints'' a.s.; this
is the case for transient random walks in $\Z^d$: see James and Peres (1996).
However, there exist transient chains without cutpoints (see James (1996)).

\medbreak \noindent {\bf Acknowledgement.} We are indebted to Itai
Benjamini, Fran\c{c}ois Ledrappier, Robin Pemantle, Tom Salisbury,
 and Jeff Steif for useful discussions. We thank the referee
for helpful comments and corrections.
\enspace

\beginreferences

Benjamini, I., Lyons, R., Peres, Y., \and Schramm, O. (2001) Uniform
spanning forests, {\it Ann.\ Probab.} {\bf 29}, 1--65.

Burkholder, D.L., Davis, B.J. \and Gundy, R.F. (1972)
Integral inequalities for convex functions of operators on martingales,
{\it Proceedings of the Sixth Berkeley Symposium on Mathematical Statistics
and Probability (Univ. California, Berkeley, Calif., 1970/1971), Vol. II:
Probability theory}, {223--240},
{Univ. California Press}, {Berkeley, Calif.}

Erd{\H{o}}s, P. \and Taylor, S.J. (1960)
{Some intersection properties of random walk paths},
{\it Acta Math. Acad. Sci. Hungar.}
{\bf 11},
{231--248}.

Fitzsimmons, P.J.  \and Salisbury, T. (1989)
Capacity and energy for multiparameter Markov processes.
{\it Ann. Inst. H. Poincar\'e Probab. Statist.}
{\bf 25}, 325--350.

H\"aggstr\"om, O. (1995) Random-cluster measures and uniform spanning
trees, {\it Stoch.\ Proc.\ Appl.} {\bf 59}, 267--275.

James, N. (1996)
{\it Ph.D. Thesis}. University of California, Berkeley.

James, N. \and Peres, Y. (1996)
Cutpoints and exchangeable events for random walks.
{\it Teor. Veroyatnost. i Primenen.} {\bf 41}, 854--868.
Reproduced in {\it Theory Probab. Appl.} {\bf 41} (1997), 666--677.

Kahane, J.-P. (1985) {\it Some Random Series of Functions.} Second
edition,  Cambridge University Press, Cambridge.

Lawler, G. (1991)
{\it Intersections of Random Walks.} Birkh\"{a}user, Boston.

Lawler, G. (1998) Loop-erased walks intersect infinitely often in
four dimensions, {\it Electronic Communications in Probability}
{\bf 3}, 35--42.

Le Gall, J.-F. \and Rosen, J. (1991)
The range of stable random walks,
{\it Ann.\ Probab.} {\bf 19}, 650--705.

Pemantle, R. (1991) Choosing a spanning tree for the integer lattice
uniformly, {\it Ann.\ Probab.} {\bf 19}, 1559--1574.

Salisbury, T. (1996)
Energy, and intersections of  Markov chains.
{\it Random Discrete Structures}, IMA Volume 76, Aldous, D. and Pemantle, R.
(Editors), Springer-Verlag, New York, 213--225.

Wilson, D.B. (1996) Generating random spanning trees more quickly than the
cover time, {\it Proceedings of the Twenty-eighth Annual ACM Symposium on the
Theory of Computing (Philadelphia, PA, 1996)}, {296--303}, {ACM}, {New York}.

Woess, W. (2000) {\it Random Walks on Infinite Graphs and Groups.}
 Cambridge Tracts in Mathematics, {\bf 138},
Cambridge University Press, Cambridge.

\endreferences

\filbreak
\begingroup
\eightpoint\sc
\parindent=0pt\baselineskip=10pt
\def\email#1{\par\qquad {\tt #1} \smallskip}
\def\emailwww#1#2{\par\qquad {\tt #1}\par\qquad {\tt #2}\smallskip}

Department of Mathematics,
Indiana University,
Bloomington, IN 47405-5701, USA
\emailwww{rdlyons@indiana.edu}
{http://php.indiana.edu/\string~rdlyons/}

and

School of Mathematics,
Georgia Institute of Technology,
Atlanta, GA 30332-0160
\email{rdlyons@math.gatech.edu}

\medskip
\smallskip

Department of Statistics, University of California,
Berkeley, CA 94720-3860, USA
\emailwww{peres@stat.berkeley.edu}
{http://www.stat.berkeley.edu/\string~peres/}

\medskip
\smallskip

Mathematics Department,
The Weizmann Institute of Science,
Rehovot 76100, Israel

and

Microsoft Research,
One Microsoft Way,
Redmond, WA 98052, USA
\emailwww{schramm@microsoft.com}
{http://research.microsoft.com/\string~schramm/}

\endgroup

\bye